\documentclass[a4paper,reqno]{amsart}

\usepackage[T1]{fontenc}
\usepackage[british]{babel}
\usepackage{amsfonts}
\usepackage{amssymb}
\usepackage{amsthm}
\usepackage{amsmath}
\usepackage[all]{xy}

\theoremstyle{plain}
\newtheorem{theorem}[subsection]{Theorem}
\newtheorem{proposition}[subsection]{Proposition}
\newtheorem{corollary}[subsection]{Corollary}

\theoremstyle{definition}

\newtheorem{definition}[subsection]{Definition}

\newtheorem{example}[subsection]{Example}
\newtheorem{examples}[subsection]{Examples}

\DeclareMathAlphabet{\mathitbf}{OML}{cmm}{b}{it}

\newcommand{\comp}{\raisebox{0.2mm}{\ensuremath{\scriptstyle{\circ}}}}
\newcommand{\defn}{\textbf}

\newcommand{\coker}{\ensuremath{\mathsf{coker\,}}}
\renewcommand{\ker}{\ensuremath{\mathsf{ker\,}}}
\newcommand{\A}{\ensuremath{\mathcal{A}}}
\newcommand{\Ab}{\ensuremath{\mathsf{Ab}}}
\newcommand{\Loop}{\ensuremath{\mathsf{Loop}}}
\newcommand{\gp}{\ensuremath{\mathsf{gp}}}
\newcommand{\PXM}{\ensuremath{\mathsf{PXMod}}}
\newcommand{\XMod}{\ensuremath{\mathsf{XMod}}}

\newcommand{\Sol}{\ensuremath{\mathsf{Sol}}}
\newcommand{\Nil}{\ensuremath{\mathsf{Nil}}}

\newcommand{\B}{\ensuremath{\mathcal{B}}}
\newcommand{\Set}{\ensuremath{\mathsf{Set}}}
\newcommand{\Gp}{\ensuremath{\mathsf{Gp}}}

\newcommand{\noproof}{\hfill \qed}
\renewcommand{\a}{\mathitbf{a}}
\renewcommand{\b}{\mathitbf{b}}

\newcommand{\x}{\mathitbf{x}}
\newcommand{\y}{\mathitbf{y}}
\newcommand{\m}{\mathitbf{m}}
\newcommand{\n}{\mathitbf{n}}
\newcommand{\p}{\mathitbf{p}}
\newcommand{\one}{\mathitbf{1}}
\newcommand{\CExt}{\ensuremath{\mathsf{CExt}}}
\newcommand{\Ext}{\ensuremath{\mathsf{Ext}}}
\newcommand{\Arr}{\ensuremath{\mathsf{Arr}}}

\newdir{>}{{}*:(1,-.2)@^{>}*:(1,+.2)@_{>}}
\newdir{<}{{}*:(1,+.2)@^{<}*:(1,-.2)@_{<}}

\begin{document}

\email{teveraer@vub.ac.be}
\author{Tomas Everaert}
\address{Vakgroep Wiskunde, Vrije Universiteit Brussel, Pleinlaan 2, 1050 Brussel, Belgium}

\email{tim.vanderlinden@uclouvain.be}
\author{Tim Van~der Linden}
\address{Centro de Matem\'atica da Universidade de Coimbra, 3001--454 Coimbra, Portugal. Institut de recherche en math\'ematique et physique, Universit\'e catholique de Louvain,
Chemin du Cyclotron~2, 1348 Louvain-la-Neuve, Belgium}



\thanks{The first author's research was supported by Fonds voor Wetenschappelijk Onderzoek (FWO-Vlaanderen). The second author's research was supported by Centro de Matem\'atica da Universidade de Coimbra and by Funda\c c\~ao para a Ci\^encia e a Tecnologia (under grant number SFRH/BPD/38797/2007). He works as \emph{charg\'e de recherches} for Fonds de la Recherche Scientifique--FNRS}

\thanks{Published as T.~Everaert and T.~Van~der Linden, \emph{Galois theory and commutators}, Algebra Universalis \textbf{65} (2011), no.\ 2, 161--177. This final publication is available at \texttt{springerlink.com}}

\title{Galois theory and commutators}

\begin{abstract}
We prove that the relative commutator with respect to a subvariety of a variety of $\Omega$-groups introduced by the first author can be described in terms of categorical Galois theory. This extends the known correspondence between the Fr\"ohlich--Lue and the Janelidze--Kelly notions of central extension. As an example outside the context of $\Omega$-groups we study the reflection of the category of loops to the category of groups where we obtain an interpretation of the associator as a relative commutator.
\end{abstract}

\maketitle

\section*{Introduction}

This article concerns the connection between universal algebra and categorical algebra which arises when the concept of \emph{relative commutator} introduced by Everaert in~\cite{EveraertCommutator} is analysed from a Galois-theoretic point of view. It may be seen as a continuation of Janelidze and Kelly's work~\cite{Janelidze-Kelly} on a general theory of \emph{central extensions}, which gives a categorical interpretation of the relative notion of central extension introduced by Fr\"ohlich~\cite{Froehlich} and Lue~\cite{Lue} in the universal-algebraic context of Higgins' $\Omega$-groups (\cite{Higgins}; see Section~\ref{Section-Omega-Groups}). We shall explain how, in a formally precise way, the relative commutator studied in~\cite{EveraertCommutator} corresponds to a two-dimensional version of those relative central extensions.

\subsection*{Relative central extensions in varieties of $\Omega$-groups}

The definition of Fr\"ohlich and Lue involves a variety of $\Omega$-groups $\A$ and a chosen subvariety $\B$ of~$\A$. Let~$I\colon{\A\to \B}$ be the left adjoint to the inclusion functor, and $\eta$ the unit of the adjunction: for $A$ in $\A$, its reflection $IA$ is $A/\!\!\sim$, where $\sim$ is the smallest congruence on~$A$ under which the quotient algebra is in $\B$; and $\eta_{A}$ is the canonical homomorphism ${A\to IA}$. Now $[-]_{\B}\colon{\A\to \A}$ denotes the \defn{variety subfunctor} associated to~$I$, which maps an object $A$ of $\A$ to the object $[A]_{\B}$ defined via the short exact sequence
\[
\xymatrix{0 \ar[r] & [A]_{\B} \ar[r]^-{\kappa_{A}} & A \ar[r]^-{\eta_{A}} & IA \ar[r] & 0,}
\]
and a morphism $a\colon{A'\to A}$ to the induced morphism $[a]_{\B}\colon{[A']_{\B}\to [A]_{\B}}$. For instance, when $\A$ is the variety $\Gp$ of all groups and $\B$ is the subvariety $\Ab$ of abelian groups, $IA$ is the abelianisation $A/[A,A]$ of $A$, and $[A]_{\Ab}$ is the commutator subgroup $[A,A]$. Of course, $A$ is in $\Ab$ if and only if $[A]_{\Ab}=[A,A]=0$.

An \defn{extension} $f\colon{A\to B}$ in $\A$ being a surjection, it is \defn{central with respect to $\B$} if and only if for any parallel pair of morphisms $a_{0}$, $a_{1}\colon{A'\to A}$, the condition $f\comp a_{0}=f\comp a_{1}$ implies that $[a_{0}]_{\B}=[a_{1}]_{\B}$. For instance, in the case of groups vs.\ abelian groups, a surjective group homomorphism $f$ is central with respect to $\Ab$ exactly when it is central in the classical sense, i.e., the commutator $[K,A]$ of $A$ with the kernel $K$ of $f$ is trivial.

\subsection*{The Galois-theoretic approach}

Janelidze and Kelly understood how this relative concept of central extension may be described in terms of \emph{categorical Galois theory}. Introduced by Janelidze (in~\cite{Janelidze:Pure}; see also~\cite{Borceux-Janelidze}) this general approach to Galois theory not only captures, e.g., the case of separable field extensions, but similar concepts in other parts of mathematics as well---indeed also Fr\"ohlich and Lue's relative central extensions, as explained in the article~\cite{Janelidze-Kelly}. In the context where we shall need it, their definition of central extension involves a semi-abelian category $\A$ (in the sense of~\cite{Janelidze-Marki-Tholen}: pointed, Barr exact and Bourn protomodular with binary sums) and a \defn{Birkhoff subcategory} $\B$ of~$\A$: full and reflective in~$\A$, and closed in $\A$ under formation of subobjects and regular quotients. Note that all varieties of $\Omega$-groups are semi-abelian and that regular quotients and regular epimorphisms in varieties of algebras are just quotients and surjective homomorphisms respectively. Moreover, here a Birkhoff subcategory  is the same thing a subvariety.

The adjunction
\begin{equation}\label{Birkhoff-Adjunction}
\xymatrix@1{{\A} \ar@<1ex>[r]^-{I} \ar@{}[r]|-{\perp} & {\B} \ar@<1ex>[l]^-{\supset}}
\end{equation}
together with the classes $|\Ext\A|$ and $|\Ext\B|$ of \defn{extensions} (i.e., regular epimorphisms) in $\A$ and in $\B$ form a so-called \emph{Galois structure}
\[
\Gamma=(\xymatrix@1{{\A} \ar@<1ex>[r]^-{I} \ar@{}[r]|-{\perp} & {\B} \ar@<1ex>[l]^-{\supset}},|\Ext\A|,|\Ext\B|).
\]
An extension $f\colon{A\to B}$ in $\A$ is \defn{trivial (with respect to $\B$)} or \defn{$\B$-trivial} when the induced commutative square
\[
\vcenter{\xymatrix{A \ar[r]^-{f} \ar[d]_-{\eta_{A}} & B \ar[d]^-{\eta_{B}}\\
IA \ar[r]_-{If} & IB}}
\]
is a pullback. And $f$ is \defn{central (with respect to $\B$)} or \defn{$\B$-central} when there exists an extension $g\colon{C\to B}$ such that the pullback $g^{*}f$ of $f$ along $g$ is $\B$-trivial.

In the present context, if $f$ is central then one necessarily has that either one of the projections $f_{0}$, $f_{1}$ in the kernel pair $(R[f],f_{0},f_{1})$ of $f$ is $\B$-trivial (i.e., $f$ is \emph{normal} with respect to $\B$)~\cite[Theorem~4.8]{Janelidze-Kelly}. When $\A$ is a variety of $\Omega$-groups, an extension is $\B$-central in the Galois-theoretic sense if and only if it is $\B$-central in the Fr\"ohlich--Lue sense~\cite[Theorem~5.2]{Janelidze-Kelly}. More precisely, the definition of central extensions as we presented it above for $\Omega$-groups is equally valid in semi-abelian categories~\cite[Theorem~2.1]{Bourn-Gran}.

\subsection*{Connections with homological algebra}

There are close connections between the Janelidze--Kelly theory of central extensions and some recent developments in homological algebra, which are worth exploring before we go deeper into the link with commutator theory. Already in the work of Fr\"ohlich, Lue and Furtado-Coelho~\cite{Froehlich, Lue, Furtado-Coelho} in the varietal context, the relation between the derived functors of the reflector $I\colon{\A\to\B}$ and the notion of $\B$-central extension is emphasised. This relation is particularly explicit in the Hopf formula
\[
H_{2}(B,\B)\cong\frac{K\cap [A]_{\B}}{[K,A]_{\B}}
\]
which gives an interpretation of a derived functor of $I$ (the left hand side of the equation) in terms of commutators (the equation's right hand side). Here the short exact sequence
\begin{equation}\label{Presentation}
\xymatrix{0 \ar[r] & K \ar[r]^-{k} & A \ar[r]^-{f} & B \ar[r] & 0}
\end{equation}
is a \defn{presentation} of $B$, i.e., $A$ is projective, and the commutator $[K,A]_{\B}$ is the smallest ideal $J$ of $A$ such that the induced map ${A/J\to B}$ is central. Similar ideas are known to work in the context of semi-abelian categories~\cite{EG2, EverVdL1, EverVdL2}.

The relative theory of central extensions may be used to prove higher-dimensional versions of the Hopf formula, which express the higher homology groups $H_{n}(B,\B)$ in terms of \emph{higher-dimensional central extensions}~\cite{EGVdL, EverHopf}. Just like the concept of central extension which is defined with respect to an adjunction~\eqref{Birkhoff-Adjunction}, one may consider double central extensions which are defined with respect to the reflection of extensions to central extensions---the adjunction
\[
\xymatrix@1{{\Ext\A} \ar@<1ex>[r]^-{I_{1}} \ar@{}[r]|-{\perp} & {\CExt_{\B}\A} \ar@<1ex>[l]^-{\supset}}
\]
where $\Ext\A$ denotes the category of extensions and commutative squares between them, and $\CExt_{\B}\A$ its full subcategory determined by those extensions which are $\B$-central. Together with well-chosen classes of \emph{double extensions}, this adjunction forms a Galois structure $\Gamma_{1}$, and Galois theory provides us with a notion of \emph{relative double central extension} with respect to $\B$. This construction may be repeated ad infinitum, so that notions of \emph{$n$-fold central extension} are obtained, but for the present purposes the second step is sufficient. Thus double central extensions, first introduced by Janelidze for groups~\cite{Janelidze:Double}, appear naturally in the study of (co)homology~\cite{EverHopf, EGVdL, EverVdL3, GVdL2, RVdL}, and turn out to be precisely what we need to understand how the relative commutator works.

\subsection*{The relative commutator}

Given a variety of $\Omega$-groups $\A$ and a subvariety $\B$ of $\A$, the objects of $\B$ and the $\B$-central extensions are often defined in terms of some kind of commutator. (Such as, in the case of groups, the classical commutator $[-,-]$, which is used to characterise abelian groups and central extensions.) One may take the opposite point of view, and ask whether the subvariety $\B$ (and the $\B$-central extensions) determine a notion of commutator. In his paper~\cite{EveraertCommutator}, the first author does exactly this: he introduces a relative commutator $[-,-]_{\B}$ with respect to $\B$ which is such that an object $A$ of $\A$ is in~$\B$ if and only if $[A,A]_{\B}=0$. Moreover, an extension $f\colon{A\to B}$ in $\A$ with kernel $K$ is $\B$-central if and only if $[K,A]_{\B}=0$. For example, if $\A$ is the variety $\PXM$ of precrossed modules (which is indeed a variety---see, e.g.,~\cite{LR}) and $\B$ is the subvariety $\XMod$ of crossed modules then~$[-,-]_{\B}$ is the \emph{Peiffer commutator} $\langle-,-\rangle$.

In view of the above ideas, a relation between this relative commutator and the concept of higher central extension was to be expected. Morally, the fact that we need \emph{double} central extensions here is simply a consequence of the commutator's having two arguments. Indeed, the objects of the Birkhoff subcategory $\B$ of~$\A$ are those objects $A$ for which the commutator $[A,A]_{\B}$ of~$A$ with itself is zero: to characterise the zero-dimensional $\B$-central extensions of $\A$, the commutator has to take no non-trivial arguments. To characterise a one-dimensional $\B$-central extension $f\colon{A\to B}$ with kernel $K$ as an extension~$f$ such that $[K,A]_{\B}$ is zero, one non-trivial argument of the commutator is needed; so it is natural to expect that a commutator with two non-trivial arguments corresponds to the two-dimensional~$\B$-central extensions of $\A$.

\subsection*{Structure of the text}

In the first section we recall the definition of the relative commutator introduced in~\cite{EveraertCommutator}, as well as some basic notation and examples. In Section~\ref{Section-Galois} we sketch the needed categorical background: semi-abelian categories, categorical Galois theory and the concept of double central extension. We characterise double central extensions in terms of double equivalence relations and the zero-dimensional commutator (Proposition~\ref{Proposition-Characterisation-Double-Extensions}). Section~\ref{Section-Varieties} contains the main result of the article, Theorem~\ref{Theorem-Characterisation}, which gives an interpretation of the relative commutator in terms of double central extensions. Finally, in Section~\ref{Section-Loops}, we study a non-classical example: the relative commutator of loops with respect to groups. The category $\Loop$ of loops and loop homomorphisms does not form a variety of $\Omega$-groups, hence lies beyond the scope of the theory introduced in~\cite{EveraertCommutator}. Nevertheless the concept of relative commutator arises naturally here when the reflection to the category $\Gp$ of groups is considered (Theorem~\ref{Theorem-Loops}).

\section{The relative commutator}\label{Section-Omega-Groups}

A \defn{variety of $\Omega$-groups}~\cite{Higgins} is a pointed variety of algebras (it has exactly one constant) that has amongst its operations and identities those of the variety of groups. It is well known that any such variety is semi-abelian~\cite{Janelidze-Marki-Tholen}. Examples include all varieties of groups, rings, modules and all kinds of algebras over rings, precrossed and crossed modules, and many others. 

We shall denote finite ordered sets $(x_{1}, x_{2}, \dots, x_{r})$, $(a_{1}, a_{2}, \dots, a_{s})$, \dots by $\x$, $\a$, etc. Instead of $(x_{1} y_{1} , x_{2} y_{2}, \dots, x_{r} y_{r})$ we shall write $\x\y$. Also $w(x_{1} , x_{2} , \dots, x_{r})$ and $w(a_{1} , a_{2} ,\dots, a_{s})$ become $w(\x)$ and $w(\a)$, for terms (words) $w$. If $1$ is the unit of a group operation we shall write $\one$ instead of $(1,\dots,1)$. Also, $a_{1} , a_{2} ,\dots, a_{s} \in A$ will be abbreviated to $\a \in A$.

In this context, a normal subobject $N$ of an object $A$ is usually called an \defn{ideal}. The meet of ideals is their intersection, and also the join is easy to compute: if $M$ and $N$ are ideals of $A$ then $M\vee N$ is the (internal) product
\[
M\cdot N=\{mn\mid m\in M,n\in N\}
\]
with the $\Omega$-group structure induced by $A$. A priori, $M\vee N$ contains all $w(\m\n)$ where $w$ is a term, $\m\in M$ and $\n\in N$; but $w(\m\n)$ may be written as the product of $w(\m\n)w(\one\n)^{-1}\in M$ with $w(\one \n)\in N$, which explains why $M\vee N=M\cdot N$. (To see that $w(\m\n)w(\one\n)^{-1}$ is indeed an element of~$M$, apply the canonical map~${A\to A/M}$ to it.) 

Given an ideal $J$ of $M\cdot N$, the canonical map $M\cdot N\to (M\cdot N)/J$ is denoted by $q_J$;
we write $R_{J}$ for the kernel pair of $q_{J}$. 

Suppose that $\A$ is a variety of $\Omega$-groups and $\B$ is a subvariety of $\A$. Then $\B$ is completely determined by a set of identities of terms of the form $w(\x) = 1$. The set of all corresponding terms $w(\x)$ forms a group
\[
W=W_{\B}=\{w(\x)\mid w(\b)=1,\forall B\in \B,\forall \b\in B\}.
\]
An object $B$ belongs to $\B$ if and only if $w(\b)=1$ for all $w\in W$ and $\b\in B$ and, consequently,
\[
[A]_{\B}=\{w(\a)\mid w\in W, \a\in A\}
\]
for any $A$ in $\A$.

Now we are in a position to recall the definition of the relative commutator introduced in~\cite{EveraertCommutator}.

\begin{definition}\label{Definition-Tomas-Commutator}
Let $\A$ be a variety of $\Omega$-groups and $\B$ a subvariety of $\A$. For any object $A$ of $\A$ and any pair of ideals $M$ and $N$ of $A$, the commutator $[M,N]_{\B}$ is the ideal of $M\cdot N$ generated by the set
\[
\{w(\m\n)w(\n)^{-1}w(\m)^{-1} w(\p)\mid w\in W, \m\in M, \n\in N, \p\in M\cap N\}.
\]
\end{definition}

\begin{examples}
We already mentioned that in the case of groups vs.\ abelian groups, $[M,N]_{\Ab}$ is the classical commutator $[M,N]$. It is shown in~\cite{EveraertCommutator} that, more generally, for any variety of $\Omega$-groups $\A$, the commutator $[M,N]_{\Ab\A}$ relative to the subvariety~$\Ab\A$ of abelian $\Omega$-groups in $\A$ is the Higgins commutator from~\cite{Higgins}. Proposition~2.3 in~\cite{EveraertCommutator} states that the commutator $[M,N]_{\XMod}$ of two precrossed submodules~$M$ and~$N$ of a precrossed module $A$ is precisely the Peiffer commutator~$\langle M,N\rangle$ of $M$~and~$N$. For any $k\geq 2$, a description of the commutator $[M,N]_{\Sol_k}$ of groups~vs.\ solvable groups of class at most $k$, and of $[M,N]_{\Nil_k}$, the commutator~of groups vs.~nilpotent groups of class at most $k$ can be found in~\cite{EGAM}.
\end{examples}

\section{Double central extensions}\label{Section-Galois}

In this section we sketch the categorical and Galois-theoretic background needed for the definition of double central extension, and we explain how those double central extensions can be characterised in terms of zero-di\-men\-sion\-al commutators and double equivalence relations (Proposition~\ref{Proposition-Characterisation-Double-Extensions}). More on these subjects may, for instance, be found in the articles~\cite{EverHopf, EGVdL, RVdL}.

\subsection{Semi-abelian categories}
Recall that a \defn{regular epimorphism} is a coequaliser of some pair of arrows. A category is \defn{regular} when it is finitely complete with coequalisers of kernel pairs and with pullback-stable regular epimorphisms. In a regular category, any morphism may be factored as a regular epimorphism followed by a monomorphism, and this \defn{image factorisation} is unique up to isomorphism. A category is \defn{Barr exact} when it is regular and such that any internal equivalence relation is \defn{effective}, i.e., it is the kernel pair of its coequaliser. Throughout this article, given a morphism $f\colon{A\to B}$, the pullback $A\times_{B}A$ (the kernel pair of $f$) is denoted by $R[f]$ and the pullback projections by~$f_{0}$ and~$f_{1}$; so that, in the varietal case, $R[f] = \{(a,a') \in A\times A\mid f(a) = f(a')\}$.

When a category is pointed and regular, \defn{Bourn protomodularity} can be defined via the \defn{regular Short Five Lemma}~\cite{Bourn1991}: for any commutative diagram
\[
\vcenter{\xymatrix{K[f'] \ar[r]^-{\ker f'} \ar[d]_-k & A' \ar[r]^-{f'} \ar[d]_-a & B' \ar[d]^-b \\ K[f] \ar[r]_-{\ker f} & A \ar[r]_-{f} & B}}
\]
such that $f$ and $f'$ are regular epimorphisms, $k$ and $b$ being isomorphisms implies that $a$ is an isomorphism. (Here, as throughout the paper, $(K[f],\ker f)$ will denote the kernel of a morpism $f$.) A \defn{semi-abelian} category is a pointed, Barr exact and Bourn protomodular category with binary coproducts~\cite{Janelidze-Marki-Tholen}. 

Since a regular epimorphism in a semi-abelian category is always the cokernel of its kernel, the following notion of (short) exact sequence is appropriate. A composable pair of morphisms
\[
\xymatrix{K \ar[r]^-{k} & A \ar[r]^-{f} & B}
\]
is \defn{exact} when the monomorphism in the image factorisation of $k$ is the kernel of~$f$. A sequence of composable morphisms is \defn{exact} when each pair of consecutive morphisms in the sequence is exact. A \defn{short exact sequence} is an exact sequence of the form
\[
\xymatrix{0 \ar[r] & K \ar[r]^-{k} & A \ar[r]^-{f} & B \ar[r] & 0;}
\]
this means that $k=\ker f$ and $f=\coker k$.

\subsection{Double extensions}\label{Subsection-Double-Extensions}
Let $\A$ be a semi-abelian category. Recall that an \defn{extension} in $\A$ is a regular epimorphism. The category $\Ext\A$ has for its objects the extensions in $\A$, and for its arrows the commutative squares between extensions. Since $\Ext\A$ need no longer be semi-abelian, we usually make all constructions involving exact sequences etc.\ in the semi-abelian category $\Arr\A$ of all arrows in $\A$ which contains $\Ext\A$ as a full subcategory.

A \defn{double extension} is a commutative square
\begin{equation}\label{Diagram-Square}
\vcenter{\xymatrix{X \ar[r]^-{c} \ar[d]_-{d} & C \ar[d]^-{g}\\
D \ar[r]_-{f} & Z}}
\end{equation}
such that all its maps and the comparison map $(d,c)\colon{X\to D\times_Z C}$ to the pullback of $f$ with $g$ are regular epimorphisms. The category of double extensions in $\A$ and commutative cubes them between will be denoted by~$\Ext^{2}\!\A$. The basic categorical properties of higher-dimensional extensions are explored in~\cite{EGVdL} and~\cite{EverHopf}.

\subsection{Double relations}\label{Subsection-Double-Relations}
Given two (internal) equivalence relations $R$ and $S$ on an object $X$, a \defn{double equivalence relation $C$ on $R$ and $S$} is an equivalence relation $C$ on $S$ of which the ``object part'' is $R$, as in the next diagram. That is, each of the four pairs of parallel morphisms on this diagram represents an equivalence relation, and these relations are compatible in an obvious sense.
\[
\vcenter{\xymatrix{C \ar@<.5ex>[r] \ar@<-.5ex>[r] \ar@<.5ex>[d] \ar@<-.5ex>[d] & S \ar@<.5ex>[d] \ar@<-.5ex>[d]\\
R \ar@<.5ex>[r] \ar@<-.5ex>[r] & X}}
\]
For instance, $R \square S$ denotes the largest double equivalence relation on $R$ and~$S$. In the special case of a variety of algebras it consists of all quadruples $(x, y, z, t)$ in $X^{4}$ in the configuration
\[
\left(\begin{matrix}
 x & y \\
 z & t
\end{matrix} \right),
\]
i.e., where $(x,z)$, $(y,t)\in R$ and $(x,y)$, $(z,t)\in S$. We shall be especially interested in the special case where $C$ is induced by a double extension~\eqref{Diagram-Square} as follows: $R=R[c]$ is the kernel pair of $c$, $S=R[d]$ is the kernel pair of $d$ and $C=R[c]\square R[d]$. It is easily seen that then the rows and columns of the induced diagram
\begin{equation}\label{Blokske}
\vcenter{\xymatrix{R[c]\square R[d] \ar@<.5ex>[r]^-{p_1} \ar@<-.5ex>[r]_-{p_0} \ar@<.5ex>[d]^-{r_1} \ar@<-.5ex>[d]_-{r_0} & R[d] \ar[r]^-{p} \ar@<.5ex>[d]^-{d_1} \ar@<-.5ex>[d]_-{d_0} & R[g] \ar@<.5ex>[d]^-{g_{1}} \ar@<-.5ex>[d]_-{g_{0}} \\
R[c] \ar@<.5ex>[r]^-{c_1} \ar@<-.5ex>[r]_-{c_0} \ar[d]_-{r} & X \ar[r]^-{c} \ar[d]_-{d} & C \ar[d]^-{g}\\
R[f] \ar@<.5ex>[r]^-{f_{1}} \ar@<-.5ex>[r]_-{f_{0}} & D \ar[r]_-{f} & Z}}
\end{equation}
are exact forks, i.e., consist of effective equivalence relations with their coequalisers.

\subsection{Birkhoff subcategories}\label{Subsection-Birkhoff} 
Given a semi-abelian category $\A$, a \defn{Birkhoff subcategory} $\B$ of~$\A$ is full and reflective in~$\A$, and closed in $\A$ under formation of subobjects and regular quotients. 

\begin{example}\label{Example-Varieties}
Recall from~\cite{Bourn-Janelidze} that a variety of algebras $\A$ is semi-abelian if and only if it has a unique constant $1$ and, for some natural number $n\geq 1$, $n$ binary terms $t_{i}$ and an $(n+1)$-ary term $t$ such that $t_{i}(x,x)=1$ and
\[
t(t_{1}(x,y),t_{2}(x,y),\dots,t_{n}(x,y),y)=x.
\]
This is the case, precisely when the variety is pointed and \emph{BIT speciale} in the sense of~\cite{Ursini2} or \emph{classically ideal determined} in the sense of~\cite{Ursini3}.

A Birkhoff subcategory $\B$ of $\A$ is the same thing as a subvariety. Since $x=y$ in $\A$ if and only if $t_{i}(x,y)=1$ for all $i$, the subvariety $\B$ is completely determined by a set of identities of terms of the form $w(\x)=1$, as in the case of varieties of $\Omega$-groups.
\end{example}

\subsection{The centralisation functor}\label{Subsection-Centralisation-Functor}

Let $\B$ be a Birkhoff subcategory of $\A$. The full subcategory of $\Ext\A$ determined by the $\B$-central extensions is denoted $\CExt_{\B}\A$. The inclusion $\CExt_{\B}\A\subset\Ext\A$ has a left adjoint, the \defn{centralisation functor}, which is denoted $I_{1}\colon{\Ext\A\to\CExt_{\B}\A}$. It may be described in terms of one-dimen\-sional commutators as follows.

Let $f\colon{A\to B}$ be an extension with kernel $K$. Note that $f$ is $\B$-central if and only if for the kernel pair $(R[f],f_{0},f_{1})$ of $f$, the morphisms
\[
[f_{0}]_{\B},[f_{1}]_{\B}\colon{[R[f]]_{\B}\to [A]_{\B}}
\]
induced by the two projections are equal: $[f_{0}]_{\B}=[f_{1}]_{\B}$. Indeed, if the latter condition holds and $f\comp a_{0}=f\comp a_{1}$ for some $a_{0}$, $a_{1}\colon{A'\to A}$, then there exists a map $a\colon{A'\to R[f]}$ such that $f_{0}\comp a=a_{0}$ and $f_{1}\comp a=a_{1}$, which implies $[a_{0}]_{\B}=[a_{1}]_{\B}$. The converse is obvious.

Since $[f_{0}]_{\B}$ and $[f_{1}]_{\B}$ are jointly monomorphic and have a common splitting, $[f_{0}]_{\B}$ is equal to $[f_{1}]_{\B}$ exactly when either one of these maps is an isomorphism. Hence the kernel $[K,A]_{\B}$ of $[f_{0}]_{\B}$ measures how far $f$ is from being central: $f$~is $\B$-central if and only if $[K,A]_{\B}$ is zero. This one-dimensional commutator $[K,A]_{\B}$ may be considered as a normal subobject of $A$ via the composite $\kappa_{A}\comp [f_{1}]_{\B}\comp \ker [f_{0}]_{\B}\colon{[K,A]_{\B}\to A}$, as displayed in the following diagram.
\begin{equation}\label{Commutator}
\vcenter{\xymatrix{& [K,A]_{\B} \ar[d]_-{\ker [f_{0}]_{\B}} \ar@{.>}[ddr] \\
0 \ar[r] & [R[f]]_{\B} \ar[r]^-{\kappa_{R[f]}} \ar@<-.5ex>[d]_-{[f_{0}]_{\B}} \ar@<.5ex>[d]^-{[f_{1}]_{\B}} & R[f] \ar@<-.5ex>[d]_-{f_{0}} \ar@<.5ex>[d]^-{f_{1}} \ar[r]^-{\eta_{R[f]}} & I R[f] \ar@<-.5ex>[d]_-{I f_{0}} \ar@<.5ex>[d]^-{I f_{1}} \ar[r] & 0\\
0 \ar[r] & [A]_{\B} \ar[r]_-{\kappa_{A}} & A \ar[r]_-{\eta_{A}} & IA \ar[r] & 0}}
\end{equation}

\begin{examples}
Here are two examples, taken from~\cite{EGVdL} and~\cite{EveraertCommutator}. If $f\colon{A\to B}$ is an extension of precrossed modules with kernel $K$ then the commutator $[K,A]_{\XMod}$ relative to the subvariety of crossed modules is the Peiffer commutator $\langle K,A\rangle$. If $f$ is an extension of groups then the commutator $[K,A]_{\Sol_{2}}$ relative to the subvariety $\Sol_{2}$ of all groups which are solvable of class at most  $2$ is $[[K,A],[A,A]]$.
\end{examples}

Given any extension $f\colon{A\to B}$ with kernel $K$, its \defn{centralisation} $I _{1}f$ is now obtained through the diagram with exact rows
\[
\xymatrix{0 \ar[r] & [K,A]_{\B} \ar[r] \ar[d] & A \ar[r]^-{\rho^{1}_{f}} \ar[d]_-{f} & \tfrac{A}{[K,A]_{\B}} \ar[r] \ar[d]^-{I _{1}f} & 0\\
& 0 \ar[r] & B \ar@{=}[r] & B \ar[r] & 0.}
\]
Considering this diagram as a short exact sequence
\[
\xymatrix{0 \ar[r] & K[\eta^{1}_{f}] \ar[r]^-{\kappa^{1}_{f}} & f \ar[r]^-{\eta^{1}_{f}} & I _{1}f \ar[r] & 0}
\]
in the semi-abelian category of arrows $\Arr\A$ we obtain a description of the unit $\eta^{1}$ of the adjunction and its kernel~$\kappa^{1}$.

Thus the Galois structure $\Gamma$ induces a new Galois structure
\[
\Gamma_{1}=(\xymatrix@1{{\Ext\A} \ar@<1ex>[r]^-{I _{1}} \ar@{}[r]|-{\perp} & {\CExt_{\B}\A} \ar@<1ex>[l]^-{\supset}},|\Ext^{2}\!\A|,|\Ext\CExt_{\B}\A|)
\]
where $\Ext\CExt_{\B}\A$ consists of all double extensions which lie in $\CExt_{\B}\A$.

\subsection{Double central extensions}\label{Subsection-Double-Central-Extensions}

By definition~\cite{EGVdL}, a double extension is a \defn{double central extension} when it is a covering~\cite{Janelidze:Pure} with respect to the Galois structure $\Gamma_{1}$. This means (cf.\ the diagrams~\eqref{Blokske} and~\eqref{Commutator}) that the double extension~\eqref{Diagram-Square}, considered as a map $(c,f)\colon{d\to g}$ in the category $\Ext \A$, is central if and only if the first projection
\[
\vcenter{\xymatrix{R[c] \ar[r]^-{c_{0}} \ar[d]_-{r} & X \ar[d]^-{d} \\
R[f] \ar[r]_-{f_{0}} & D}}
\qquad\qquad
\vcenter{\xymatrix{R[c] \ar[r]^-{c_{0}} \ar[d]_-{\rho^{1}_{r}} & X \ar[d]^-{\rho^{1}_{d}} \\
\tfrac{R[c]}{[K[r],R[c]]_{\B}} \ar[r] & \tfrac{X}{[K[d],X]_{\B}}}}
\]
of its kernel pair---the left hand side square---is a trivial extension with respect to~$\Gamma_{1}$. (Alternatively, one could use the square of second projections.) This means that the comparison map to its reflection into $\CExt_{\B}\A$---the right hand side square---is a pullback. For this to happen, the natural map
\begin{equation}\label{Comparison-Map}
{[K[r],R[c]]_{\B}\to [K[d],X]_{\B}}
\end{equation}
must be an isomorphism. This, in turn, is equivalent to the square
\[
\vcenter{\xymatrix{[R[c]\square R[d]]_{\B} \ar[r]^-{[p_0]_{\B}} \ar[d]_-{[r_0]_{\B}} & [R[d]]_{\B} \ar[d]^-{[d_0]_{\B}} \\
[R[c]]_{\B} \ar[r]_{[c_0]_{\B}} & [X]_{\B}}}
\]
being a pullback, because $[K[r],R[c]]_{\B}$ and $[K[d],X]_{\B}$ are by definition the kernels of the vertical maps in this square (compare diagrams \eqref{Blokske} and \eqref{Commutator}). Note that, equivalently, we could have chosen second projections for the vertical maps. Thus we proved the following characterisation of double central extensions, which is a relative version of the result in Section~1.8 of~\cite{RVdL}.

\begin{proposition}\label{Proposition-Characterisation-Double-Extensions}
Let $\A$ be a semi-abelian category and $\B$ a Birkhoff subcategory of $\A$. A double extension~\eqref{Diagram-Square} is central with respect to $\B$ if and only if any of the induced commutative squares in
\[
\vcenter{\xymatrix{[R[c]\square R[d]]_{\B} \ar@<-.5ex>[r] \ar@<.5ex>[r] \ar@<-.5ex>[d] \ar@<.5ex>[d] & [R[d]]_{\B} \ar@<-.5ex>[d] \ar@<.5ex>[d] \\
[R[c]]_{\B} \ar@<-.5ex>[r] \ar@<.5ex>[r] & [X]_{\B}}}
\]
is a pullback.\noproof
\end{proposition}

\section{The commutator in terms of categorical Galois theory}\label{Section-Varieties}

Now we are ready to prove our main theorem: a characterisation of the relative commutator from Definition~\ref{Definition-Tomas-Commutator} in Galois-theoretic terms.

\begin{theorem}\label{Theorem-Characterisation}
Let $\A$ be a variety of $\Omega$-groups and $\B$ a subvariety of $\A$. Given any two ideals $M$ and $N$ of an object $A$ of $\A$, the commutator $[M,N]_{\B}$ is zero if and only if the double extension
\begin{equation}\label{Double-Extension}
\vcenter{\xymatrix{M\cdot N \ar[r]^-{q_{M}} \ar[d]_-{q_{N}}  & \tfrac{M\cdot N}{M} \ar[d]\\
\tfrac{M\cdot N}{N} \ar[r] & 0}}
\end{equation}
is central.
\end{theorem}
\begin{proof}
First note that the above square is indeed a double extension: the comparison map
\[
(q_{N},q_{M})\colon{M\cdot N\to \tfrac{M\cdot N}{N}\times \tfrac{M\cdot N}{M}}
\]
to the induced pullback is a surjection, because
\begin{align*}
(q_{N}(mn),q_{M}(m'n')) &= (q_{N}(m),q_{M}(n'))\\
&=(q_{N}(mn'),q_{M}(mn'))
\end{align*}
for all $m$, $m'\in M$, $n$, $n'\in N$.

Now suppose that $[M,N]_{\B}$ is zero. By Proposition~\ref{Proposition-Characterisation-Double-Extensions}, it suffices to prove that any of the commutative squares in
\begin{equation}\label{Square-Commutator}
\vcenter{\xymatrix{[R_{M}\square R_{N}]_{\B} \ar@<.5ex>[r] \ar@<-.5ex>[r] \ar@<.5ex>[d] \ar@<-.5ex>[d] & [R_{N}]_{\B} \ar@<.5ex>[d] \ar@<-.5ex>[d]\\
[R_{M}]_{\B} \ar@<.5ex>[r] \ar@<-.5ex>[r] & [M\cdot N]_{\B}}}
\end{equation}
is a pullback. Now $R_{M}\square R_{N}$ consists of all 
\[
\left(\begin{matrix}
mn & mn'\\
m'n & m'n'p
\end{matrix} \right)
\]
where $m$, $m'\in M$, $n$, $n'\in N$ and $p\in M\cap N$. Hence $[R_{M}\square R_{N}]_{\B}$ contains all quadruples 
\[
\left(\begin{matrix}
w(\m\n) & w(\m\n')\\
w(\m'\n) & w(\m'\n'\p)
\end{matrix} \right)
\]
where $w\in W$, $\m, \m'\in M$, $\n, \n'\in N$ and $\p\in M\cap N$. We have to prove that two of those quadruples coincide as soon as three out of four of their elements do. This is the case because the assumption implies 
\begin{align*}
w(\m'\n'\p) &= w(\m')w(\n')w(\p)\\
&= w(\m')w(\n')\\
&= w(\m')w(\n)w(\n)^{-1}w(\m)^{-1}w(\m)w(\n')\\
&= w(\m'\n)w(\m\n)^{-1}w(\m\n'),
\end{align*}
so that the fourth element depends fully on the other three.

Conversely, suppose that any of the commutative squares in \eqref{Square-Commutator} is a pullback. Then, since for any $w\in W$ and $\p\in M\cap N$, both
\[
{\left(\begin{matrix}
1 & 1\\
1 & w(\p)
\end{matrix} \right)}\qquad\text{and}\qquad{\left(\begin{matrix}
1 & 1\\
1 & 1
\end{matrix} \right)}
\]
are in $[R_{M}\square R_{N}]_{\B}$, we have that $w(\p)=1$. Also, for any $w\in W$, $\m\in M$ and~${\n'\in N}$, we have that both
\[
{\left(\begin{matrix}
w(\m) & w(\m\n')\\
1 & w(\n')
\end{matrix} \right)}\qquad\text{and}\qquad{\left(\begin{matrix}
w(\m) & w(\m\n')\\
1 & w(\m)^{-1}w(\m\n')
\end{matrix} \right)}
\]
are in $[R_{M}\square R_{N}]_{\B}$, hence $w(\m\n')=w(\m)w(\n')$.
\end{proof}

Using Theorem~\ref{Theorem-Characterisation}, one can extend the relative notion of commuting subobjects so that $\A$ is allowed to be any semi-abelian category and $\B$ any Birkhoff subcategory of $\A$. Taking, for normal subobjects $M$ and $N$ of an object $A$ in~$\A$, their commutator $[M,N]_{\B}$ to be the smallest normal subobject $J$ of $M\vee N$ such that $q_J M$ and $q_J N$ commute thus provides one with a categorical notion of relative commutator, which will be studied in more detail in~\cite{EverVdLRCT}.

\section{Example: the associator of loops}\label{Section-Loops}

We now illustrate this approach with the example of loops vs.\ groups. Note that loops  do not constitute a variety of $\Omega$-groups, so that Definition~\ref{Definition-Tomas-Commutator} is not applicable. Nevertheless, the variety $\Loop$ is semi-abelian. The commutator $[-,-]_{\B}$ defined above in terms of double central extensions characterises the associator of loops when $\B$ is taken to be the subvariety $\Gp$ of groups. Since indeed loops are ``non-associative groups'' it makes sense for the reflection to $\Gp$ to induce a commutator $[M,N]_{\Gp}$ which measures how well the elements of two given normal subloops $M$ and $N$ of a loop $A$ associate with each other.

\subsection{Basic definitions and properties}

Recall that a \defn{loop} is an algebra
\[
(A,\cdot,\backslash,/,1)
\]
where the multiplication $\cdot$ and the left and right division $\backslash$ and $/$ satisfy the axioms
\begin{align*}
y &= x\cdot(x\backslash y) \qquad & y = x\backslash(x\cdot y)\\
x &= (x/y)\cdot y & x = (x\cdot y)/y
\end{align*} 
and $1$ is a unit for the multiplication, $x\cdot 1 = x = 1\cdot x$. We shall sometimes write $xy$ for $x\cdot y$. The variety of loops is denoted by $\Loop$. It is known to be semi-abelian (as mentioned for instance in~\cite{Borceux-Clementino}) and easily seen to be such using the description recalled in Example~\ref{Example-Varieties}: take $n=1$, $t(x,y)=xy$ and $t_{1}(x,y)=x/y$. 

Suppose that $M$ and $N$ are normal subloops of a loop $A$. Then the argument showing that $M\vee N=M\cdot N$ given in the case of $\Omega$-groups is still valid, as indeed~$w(\m\n)$ is the product of $w(\m\n)/w(\n)$ with $w(\n)$.

\subsection{The associator}

The \defn{associator} of three elements $x$, $y$, $z$ of a loop $A$ is the unique element $[x,y,z]$ of $A$ such that $(xy)z=[x,y,z]\cdot x(yz)$. Hence $[x,y,z]$ is equal to $((xy)z)\slash(x(yz))$. Given three normal subloops $L$, $M$ and $N$ of $A$, we write $[L,M,N]$ for the normal subloop of $L\cdot M\cdot N$ generated by the elements $[x,y,z]$, where either $(x,y,z)$ or any of its permutations is in $L\times M\times N$: it is the smallest normal subloop $J$ of $A$ such that $q_J L$, $q_J M$ and~$q_J N$ ``associate''.

The \defn{associator of $A$} is its normal subloop $[A,A,A]$. A loop $A$ is a group if and only if its associator is trivial, and the reflection $\gp A$ of any loop $A$ into $\Gp$ is given by $A/[A,A,A]$. Thus we see that $[A]_{\Gp}=[A,A,A]$.

\subsection{Characterisation of the $\Gp$-central extensions of loops}

The adjunction
\[
\xymatrix@1{{\Loop} \ar@<1ex>[r]^-{\gp} \ar@{}[r]|-{\perp} & {\Gp} \ar@<1ex>[l]^-{\supset}}
\]
induces a notion of central extension of loops, relative to the subvariety of groups. It may be characterised in terms of an associator as follows.

\begin{proposition}\label{Proposition-Central-Extensions-of-Loops}
In the variety of loops, let $f\colon{A\to B}$ be an extension with kernel~$K$. The extension $f$ is central with respect to $\Gp$ if and only if the associator $[K,A,A]$ is zero. Hence
\[
[K,A]_{\Gp}=[K,A,A].
\]
\end{proposition}
\begin{proof}
By definition, $f$ is a central extension if and only if the induced split epimorphisms of loops
\[
[f_{0}]_{\Gp},[f_{1}]_{\Gp}\colon{[R[f],R[f],R[f]]\to[A,A,A]}
\]
are equal.

Let $(k,a,a')$ be an element of $K\times A\times A$ and suppose that $[f_{0}]_{\Gp}=[f_{1}]_{\Gp}$. Then
\[
[k,a,a']=f_{0}[(k,1),(a,a),(a',a')]=f_{1}[(k,1),(a,a),(a',a')]=[1,a,a']=1.
\]
Similarly, also $[a,k,a']=1$ and $[a,a',k]=1$, which means that $[K,A,A]=0$.

Conversely, assume that $[K,A,A]=0$. First note that this implies that 
\begin{equation}\label{cancellation}
(ak)/(a'k)=a/a' \qquad \forall a,a'\in A, k\in K.
\end{equation}
Indeed, it follows from the assumption that 
\[
ak=((a/a')a')k=(a/a')(a'k).
\] 

Now, for any element $(a,a')$ of $R[f]$ we can write $a'=ak$, where $k=a\backslash a'$ is in $K$. Thus we see that $[R[f],R[f],R[f]]$ is generated by all elements of the form $[(a,ak),(a',a'k'),(a'',a''k'')]$ where $(a,a',a'')\in A^{3}$ and $(k,k',k'')\in K^{3}$. We have to show for any such generator that $[a,a',a'']=[ak,a'k',a''k'']$. We shall do this in three steps, first eliminating $k''$, then $k'$, and finally $k$.

Consider an associator $[ak,a'k',a''k'']$. Then
\[
[ak,a'k',a''k'']= ((ak\cdot a'k')(a''k''))\slash ((ak)(a'k'\cdot a''k''))  
\]
and $(ak\cdot a'k')(a''k'')=(ak\cdot a'k')a''\cdot k''$ while
\[
(ak)(a'k'\cdot a''k'')=(ak)((a'k'\cdot a'')k'')=(ak)(a'k'\cdot a'')\cdot k''
\]
so that $[ak,a'k',a''k'']=[ak,a'k',a'']$ by \eqref{cancellation}. Write $k'a''=a''\overline{k'}$ where $\overline{k'}=a''\backslash (k'a'')\in K$. Since
\[
(ak\cdot a'k')a''=((ak\cdot a') k')a'' =(ak\cdot a') (k'a'')=  (ak\cdot a')(a''\overline{k'})=(ak\cdot a') a''\cdot\overline{k'}
\]
and
\[
(ak)(a'k'\cdot a'')=(ak)(a'\cdot k'a'')=(ak)(a'\cdot a''\overline{k'})=(ak)(a'a''\cdot \overline{k'})=(ak\cdot a'a'')\cdot \overline{k'},
\]
using \eqref{cancellation} we find that $[ak,a'k',a'']=[ak,a',a'']$. Finally, 
\begin{align*}
(ak\cdot a')a''&=(a\cdot ka')a''=(a\cdot a'\overline{k})a''=(aa'\cdot \overline{k})a''\\
&=aa'\cdot \overline{k}a''=aa'\cdot a''\overline{\overline{k}}=(aa'\cdot a'')\cdot \overline{\overline{k}}
\end{align*}
and
\begin{align*}
ak\cdot a'a''&=a(k\cdot a'a'')=a(ka'\cdot a'')=a(a'\overline{k}\cdot a'')=a(a'\cdot \overline{k}a'')\\
&=a(a'\cdot a''\overline{\overline{k}})
=a(a'a''\cdot \overline{\overline{k}})=(a\cdot a'a'')\cdot \overline{\overline{k}},
\end{align*}
for some $\overline{k}, \overline{\overline{k}}\in K$, from which we infer---again using \eqref{cancellation}---that the equality
\[
[a,a',a'']=[ak,a'k',a''k'']
\]
holds. 
\end{proof}

One way to apply this result occurs when computing of the homology of loops with respect to groups: the second homology group of a loop $B$ may be written in a Hopf formula as a quotient of associators.

\begin{corollary}[Hopf formula for loops vs.\ groups]\label{Corollary-Hopf-Formula}
If $B$ is a loop and \eqref{Presentation} a projective presentation of $B$, then 
\[
H_{2}(B,\Gp)\cong\frac{K\cap [A,A,A]}{[K,A,A]},
\]
where the left hand side homology group is the comonadic homology of $B$ with coefficients in the reflector $\gp\colon{\Loop\to \Gp}$ and relative to the comonad induced by the forgetful/free adjunction to $\Set$.
\end{corollary}
\begin{proof}
This is an instance of~\cite[Theorem~6.9]{EverVdL2}; see also~\cite{EGVdL} and~\cite{EG2}.
\end{proof}

\subsection{The relative commutator is an associator}

We now use the characterisation of $\Gp$-central extensions of loops to interpret the commutator $[-,-]_{\Gp}$ in terms of associator elements. Let $A$ be a loop and let $M$ and $N$ be normal subloops of~$A$. As in the case of $\Omega$-groups, the induced commutative square \eqref{Double-Extension} is indeed a double extension of loops: the proof given with Theorem~\ref{Theorem-Characterisation} is still valid. We have to find out when this double extension is $\Gp$-central. 

\begin{theorem}\label{Theorem-Loops}
If $M$ and $N$ are normal subloops of a loop $A$ then
\[
[M,N]_{\Gp}=[M,N,M\cdot N].
\] 
\end{theorem}
\begin{proof}
We have to show that the associator $[M,N,M\cdot N]$ is zero if and only if \eqref{Double-Extension} is central. This happens when the two projections
\[
\bigl[K[q_{N}\comp (q_{M})_{0}]\cap K[q_{N}\comp (q_{M})_{1}],R_{M}\bigr]_{\Gp}\to [N,M\cdot N]_{\Gp}
\]
arising through the diagram
\[
\vcenter{\xymatrix{K[(q_{N}\comp (q_{M})_{0},q_{N}\comp (q_{M})_{1})] \ar@<.5ex>[r] \ar@<-.5ex>[r] \ar[d] & N \ar[d] \\
R_{M} \ar@<.5ex>[r]^-{(q_{M})_1} \ar@<-.5ex>[r]_-{(q_{M})_0} \ar[d]_-{(q_{N}\circ (q_{M})_{0},q_{N}\circ (q_{M})_{1})} & M\cdot N \ar[r]^-{q_{M}} \ar[d]_-{q_{N}} & \tfrac{M\cdot N}{M} \ar[d]\\
\tfrac{M\cdot N}{N}\times\tfrac{M\cdot N}{N} \ar@<.5ex>[r] \ar@<-.5ex>[r] & \tfrac{M\cdot N}{N} \ar[r] & 0}}
\]
and corresponding to \eqref{Comparison-Map} are equal to each other. Proposition~\ref{Proposition-Central-Extensions-of-Loops} tells us that they are the restrictions of the kernel pair projections $(q_{M})_{0}$ and $(q_{M})_{1}$ to maps
\[
q_{0},q_{1}\colon\bigl[K[q_{N}\comp (q_{M})_{0}]\cap K[q_{N}\comp (q_{M})_{1}],R_{M},R_{M}\bigr]\to [N,M\cdot N,M\cdot N].
\]

If $q_{0}=q_{1}$ then, for any $m$, $m'\in M$ and $n$, $n'\in N$,
\begin{align*}
[n,m,m'n']&=q_{0}[(n,n),(m,1),(m'n',n')]\\
&=q_{1}[(n,n),(m,1),(m'n',n')]\\
&=[n,1,n']=1.
\end{align*}
Similarly, $[m,n,m'n']$, $[m,m'n',n]$, $[n,m'n',m]$, $[m'n',m,n]$ and $[m'n',n,m]$ also vanish so that $[M,N,M\cdot N]=0$.

Conversely, suppose that $[M,N,M\cdot N]$ is zero. We have to prove that then the morphism $q_{0}$ is equal to $q_{1}$. Note that any element of the intersection
\[
K[q_{N}\comp (q_{M})_{0}]\cap K[q_{N}\comp (q_{M})_{1}]
\]
may be written as $(n,np)$ with $n\in N$ and $p\in M\cap N$, and any element of $R_M$ as $(nm,nmm')$ with $m, m'\in M$ and $n\in N$. Hence it will be sufficient to prove that the identities 
\begin{align*}
[n,n'm,n''m'']&=[np,n'mm',n''m''m'''],\\
[n''m'',n,n'm]&=[n''m''m''',np,n'mm']
\end{align*}
and 
\[
[n'm,n''m'',n]=[n'mm',n''m''m''',np]
\]
hold for all $p\in M\cap N$, $m$, $m'$, $m''$, $m'''\in M$ and $n$, $n'$, $n''\in N$. This in its turn simplifies to proving that each of the commutators $[nm,n'm',n'']$, $[nm,n',n''m'']$ and $[n,n'm',n''m'']$ is equal to $[n,n',n'']$, since in that case we have that
\[
 [n,n'm,n''m'']= [n,n',n'']=[np,n',n'']=[np,n'mm',n''m''m''']
\]
---and the other two identities can be derived in a similar fashion. We shall work out in detail that $[nm,n'm',n'']=[n,n',n'']$. The proofs of other two identities are quite similar and are left to the reader. 

Let thus $m$ and $m'$ be elements of $M$ and let $n$, $n'$ and $n''$ be elements of~$N$. Write $m'n''=n''\overline{m'}$,  $mn'=n'\overline{m}$ and then $\overline{m}n''=n''\overline{\overline{m}}$, where $\overline{m'}=n''\backslash (m'n'')$, $\overline{m}=n'\backslash (mn')$ and $\overline{\overline{m}}=n''\backslash (\overline{m}n'')$ all lie in $M$. Then
\[
[nm,n'm',n'']=((nm\cdot n'm')\cdot n')/(nm\cdot (n'm'\cdot n''))
\]
and
\begin{align*}
(nm\cdot n'm')n''&=((nm\cdot n')m')n''=(nm\cdot n')\cdot m'n''=(nm\cdot n')\cdot n''\overline{m'}\\
&=((nm\cdot n')n'')\overline{m'}=((n\cdot mn')n'')\overline{m'}=((n\cdot n'\overline{m})n'')\overline{m'}\\
&=((nn'\cdot \overline{m})n'')\overline{m'}=(nn'\cdot \overline{m}n'')\overline{m'}
=(nn'\cdot n''\overline{\overline{m}})\overline{m'}\\
&=((nn'\cdot n'')\overline{\overline{m}})\overline{m'}=(nn'\cdot n'')\cdot \overline{\overline{m}}\overline{m'}
\end{align*}
while
\begin{align*}
nm\cdot (n'm'\cdot n'')&=nm\cdot (n'\cdot m'n'')=nm\cdot (n'\cdot n''\overline{m'})=nm\cdot (n'n''\cdot \overline{m'})\\
&= (nm\cdot n'n'')\overline{m'}=(n(m\cdot n'n''))\overline{m'}=(n(mn'\cdot n''))\overline{m'}\\
&=(n(n'\overline{m}\cdot n''))\overline{m'}=(n(n'\cdot\overline{m} n''))\overline{m'}=(n(n'\cdot n''\overline{\overline{m}}))\overline{m'}\\
&=(n(n'n''\cdot\overline{\overline{m}}))\overline{m'}=((n\cdot n'n'')\overline{\overline{m}})\overline{m'}=(n\cdot n'n'')\cdot \overline{\overline{m}}\overline{m'}
\end{align*}
so that 
\[
(nn'\cdot n'')\cdot \overline{\overline{m}}\overline{m'}=([n,n',n'']\cdot (n\cdot n'n''))\cdot \overline{\overline{m}}\overline{m'}=[n,n',n'']\cdot ((n\cdot n'n'')\cdot \overline{\overline{m}}\overline{m'})
\]
which implies that $[nm,n'm',n'']=[n,n',n'']$.
\end{proof}



\providecommand{\bysame}{\leavevmode\hbox to3em{\hrulefill}\thinspace}

\end{document}